\documentclass{winnower}
\usepackage[T1]{fontenc}
\usepackage{ae,aecompl}
 \usepackage{tikz}
\usetikzlibrary{calc,patterns}
\usetikzlibrary{arrows,shapes,positioning}
\usetikzlibrary{decorations.markings}
\tikzstyle arrowstyle=[scale=1]
\tikzstyle directed=[postaction={decorate,decoration={markings,
    mark=at position .65 with {\arrow[arrowstyle]{stealth}}}}]
\tikzstyle reverse directed=[postaction={decorate,decoration={markings,
    mark=at position .65 with {\arrowreversed[arrowstyle]{stealth};}}}]

\usepackage{float}
\usepackage{graphicx, tabularx, multirow, psfrag, booktabs}
\usepackage[numbers]{natbib}

\begin{document}
%\includepdf[pages=1-last]{artigoTeste.pdf}

\graphicspath{figs/}
\title{\textbf{On the verification of adaptive three-dimensional multiresolution computations of the magnetohydrodynamic equations}}

\author{Anna Karina Fontes Gomes\footnote{Corresponding author.\\Email address:
anna.gomes@inpe.br }}
\affil{\small Post-graduation program in Applied Computing, National Institute for Space Research, S\~{a}o Jos\'{e} dos Campos, Brazil}
\author{Margarete Oliveira Domingues}
\affil{\small Associated Laboratory for Computing and Applied Mathematics, National Institute for Space Research, S\~{a}o Jos\'{e} dos Campos, Brazil}
\author{Odim Mendes}
\affil{\small Space Geophysics Division, National Institute for Space Research, S\~{a}o Jos\'{e} dos Campos, Brazil}
\author{Kai Schneider}
\affil{\small Institut de Math\'ematiques de Marseille, Aix-Marseille Universit\'e, CNRS, Centrale Marseille, I2M UMR 7373}

\date{}

\maketitle

\begin{abstract}
Magnetohydrodynamics is an important tool to study the dynamics of plasma Space Physics. In this context, we introduce a three-dimensional magnetohydrodynamic solver with divergence-cleaning in the adaptive multiresolution CARMEN code. The numerical scheme is based on a finite volume discretization that ensures the conservation of physical quantities. The adaptive multiresolution approach allows for automatic identification of local structures in the numerical solution and thus provides an adaptive mesh refined only in regions where the solution needs more improved resolution. We assess the three-dimensional magnetohydrodynamic CARMEN code and compare its results with the ones from the well-known FLASH code.
\end{abstract}

\section{Introduction}
\noindent The Magnetohydrodynamic (MHD) theory helps in understanding the behavior of the macroscopic space plasma scenario. 
In particular, MHD modeling describes the flow dynamics of a conducting fluid under the influence of a magnetic field and its nonlinear interaction.

The MHD model can be formulated as an initial value problem composed by four nonlinear, vector valued partial differential equations. 
In most practical cases, it does not have an analytic solution and we need numerical methods to solve it.
For instance, MHD simulations in three-dimensions can be found in \cite{otto19903d, touma2006central}  and in more recent works \cite{lee2013solution, skala20153d} using either finite difference or finite volume methods with and without adaptive mesh refinement techniques. 
Due to the high computational cost of three-dimensional ideal  MHD simulations, it is adequate to use, when possible, an adaptive mesh procedure to create a mesh which has more cells at locations where local structures are present, and there are less cells in smooth regions. 

Here, we introduce a verification of a new code for an adaptive ideal three-dimensional MHD model discretized with a finite volume method in which the domain is divided into small volumes or cells and a multiresolution technique is used to introduce an adapted mesh in order to reduce the number of operations and memory while preserving the precision of the discretization. 
Previous results were obtained for two-dimensional MHD equations in \cite{Gomes:2012:AnMuAd,domingues2013extended,Gomesetal:2015}.
To perform the present simulations we extended the adaptive multiresolution CARMEN code \cite{RSTB03} and modified the MHD equations in three-dimensions \cite{Gomes:2012:AnMuAd, Gomesetal:2015}. 
As a benchmark for our results we use the well-known astrophysical FLASH code\footnote{http://flash.uchicago.edu/} which uses a finite volume method on a regular Cartesian mesh. 
We present results for two Riemann problems in order to check our new adaptive MHD implementation.

%In the context of mesh adaptivity, we introduce the adaptive multiresolution analysis (MR) for cell-averages \cite{Harten:1996, Paper OR 2003}, which allows the mesh to be adapted to the solution of the problem by evaluating a wavelet transform to represent the data in diffent levels of resolution. Results obtained for MHD equations with MR approach in one and two dimensions were published in \cite{Gomes:2012:AnMuAd,domingues2013extended,Gomesetal:2015}.
%

The content of this paper is organized as follows: in Section~\ref{sec:mhd}, we introduce the equations of the ideal MHD model. The numerical approach to simulate the MHD equations is described in Section~\ref{sec:num}, where we also briefly describe the parabolic-hyperbolic divergence cleaning, the finite volume, and the MR methods. In Section~\ref{sec:results}, results are presented and discussed. 
Lastly, we make our conclusions and give some perspectives.

\section{The MHD model}
\label{sec:mhd}
The ideal MHD model describes the physics of an ideal conducting fluid under the influence of a magnetic field. 
The study of the MHD model and its simulation is important to better understand how an ideal plasma behaves. It enables us to perform numerical investigations of more realistic physical phenomena, which are directly related to space sciences.
%MD: reference to be included
Here we present the set of four partial differential equations that constitutes the ideal MHD model.
            \begin{subequations}
\begin{eqnarray}
	     % \label{Cont}
	    \frac{\partial\rho}{\partial t}+\nabla\cdot(\rho {\bf u})&=&0, \\
	     \label{eq:cont}
		\frac{\partial \mathcal{E}}{\partial t} +\nabla\cdot\left[ \left(\mathcal{E} + p + \frac{B^2}{2}\right){\bf u} - {\bf B}\cdot\left({\bf u}\cdot{\bf B}\right)\right] &=& 0,  \\
	      \label{eq:energy}
		\frac{\partial \rho{\bf u}}{\partial t} + \nabla\cdot\left[ \rho{\bf uu} +\left( p + \frac{1}{2} B^2 \right){\bf I} - {\bf BB} \right] &=&0, \\
	      \label{eq:mom}
		\frac{\partial{\bf B}}{\partial t} +\nabla\cdot\left( {\bf uB} - {\bf Bu}\right) &=& 0,	
\end{eqnarray}
\label{sys:ideal}
\end{subequations}
%and ${\bf J}=\nabla\times{\bf B}$ the current density. 
with the model variables: $\rho$ density, ${\bf u} = (u_x, u_y,u_z)$ the fluid velocity vector, ${\bf B}=(B_x,B_y,B_z)$ the magnetic field vector and $\mathcal{E}$ total mass density defined as
\begin{equation}
    \mathcal{E} = \frac{p}{\gamma - 1} + \frac{\rho u^2}{2} + \frac{B^2}{2},
\end{equation}
where $p$ is the pressure, and $u$ and $B$ are the magnitude of $\bf u$ and $\bf B$, respectively.

The model is written in its conservation form, \textit{i.e.}, a differential equation of the form
\begin{equation}
    \frac{\partial {\bf U}}{\partial t} + \nabla\cdot{\bf F} = {\bf 0},
    \label{eq:consLaw}
\end{equation}
with the vector of quantities $\bf U$ and the flux tensor ${\bf F}={\bf F}({\bf U})$. 
In our case, we have
\begin{equation}	
	%\begin{eqnarray}
	%&
	{\bf U} =\left(
	\begin{array}{c}
	\rho \\
	\mathcal{E}\\
	\rho{\bf u}\\{\bf B}
	\end{array}
	\right),
	%\\
	\qquad
	%\\
	%&
	 {\bf F(U)} = \left(
	\begin{array}{c}
	\rho{\bf u} \\
	\displaystyle (\mathcal{E} + p + \frac{{\bf B\cdot B}}{2}){\bf u} - {\bf B}({\bf u\cdot B})\\
	\displaystyle\rho{\bf uu} + {\bf I}(p + \frac{{\bf B \cdot B}}{2})-{\bf BB}\\
	{\bf uB-Bu}
	\end{array}
	\right).
		 \label{eq:fluxvec}
	 \end{equation}
The conservative formulation of the MHD equations describes the conservation of mass density, total energy density, momentum, and magnetic flux. 
It is important to note that the Gauss law of magnetism is also part of the model and can be understood as a physical constraint of the magnetic field. It implies that the divergence of the magnetic field is zero,  i.e.,  there are no magnetic monopoles in the solution of the MHD equations \cite{BrBa80}. 

\section{Numerical Approach}
\label{sec:num}
When it comes to numerical simulations, where the domain of the physical problem has to be discretized, the divergence constraint is not satisfied anymore due to approximation errors. 
In order to minimize these errors, it is acceptable to introduce a correction to the model so that the divergence errors are well controlled. 
There are several types of corrections and some of them are discussed here \cite{Toth:2000JCP,toth2012adaptive}. 
We have chosen a divergence cleaning approach, which does not impose zero divergence, but controls the errors such that non-physical behavior is avoided in the numerical solution of the model.%, and its implementation is simple to add in existing codes.

%\subsection{The divergence cleaning approach}
The parabolic-hyperbolic divergence cleaning was first proposed in \cite{Dedneretal:2002} and it offers both propagation to the boundary and dissipation of the divergence errors. 
The correction consists in adding a differential operator $D=\frac{1}{c_h^2}\frac{\partial}{\partial t} + \frac{1}{c_p^2}$ to the divergence constraint of the $\textbf B$ field, resulting in the Generalized Lagrange Multiplier MHD model (GLM--MHD), composed by Eqs.~(\ref{eq:cont}),~(\ref{eq:energy}),~(\ref{eq:mom}) and the additional equations
            \begin{subequations}
        	  \begin{eqnarray}
        		\displaystyle\frac{\partial{\bf B}}{\partial t} + \nabla\cdot\left({\bf uB-Bu}+\psi {\bf I}\right)&=&0, \label{eq:mag} \\
        			      \frac{\partial\psi}{\partial t}+c_h^2\nabla\cdot{\bf B}&=&-\frac{c_h^2}{c_p^2}\psi        \label{eq:div}
        	  \end{eqnarray}
        	  \label{sys:glm}
        	\end{subequations} 
where $\psi$ is a scalar-valued function, $\textbf{I}$ the identity tensor, $c_p$ and $c_h$ are the parabolic and hyperbolic constants, respectively. We can note that System~\ref{sys:glm} reduces to the ideal MHD model when $\psi=0$. 
The GLM--MHD model is used for the numerical simulation and the divergence cleaning is implemented in the flux as a source term in Equation~\ref{eq:div}. 
%The last case we have, if we consider Equation~\ref{eq:div} as
%\begin{equation}
%    \frac{\partial\psi}{\partial t}=-\frac{c_h^2}{c_p^2}\psi, 
%\end{equation}
%it is possible to obtain an analytical solution for $\psi$. In the numerical context, to compute the value of $\psi$ in time $t^{n}$, $n\in\mathbb{N}$, \textit{i. e.},
%\begin{equation}
%    \psi^n = \psi^\star\exp{-\Delta t^n\frac{c_h^2}{c_p^2}},
%\end{equation}
%it is required the time step $\Delta t^n$ and $\psi^\star$, which is the variable $\psi$ evolved in time. Basically, we update the variable with the correction term after time evolution. 
The constant $c_h$ is defined as
\begin{equation}
    c_h = \max\left[\nu\frac{\Delta h}{\Delta t},\,\max\left( |u_d|\pm c_f \right)\right],
\end{equation}
where $\Delta h = \min\,\left( \Delta x, \Delta y, \Delta z \right)$, with $\Delta x, \Delta y, \Delta z$ are the space steps in each direction, $\nu$ the Courant number, $u_d$ is the velocity of the $d$-th component and $c_f$ is the fast magnetoacoustic wave speed of the MHD model. 
The $c_p$ value is defined in terms of the parameter $\alpha=\Delta h\frac{c_h}{c_p^2}$, as discussed in \cite{mignone2010second}.

\subsection{Finite volume discretization}
The discretization of the model is performed with a finite volume method, which is based on the integral form of the conservation laws. 
The idea here is to divide the computational domain into mesh cells, assigning a cell average to the vector quantity $\overline{\textbf{U}} = \overline{\textbf{U}}_{i,j,k}$ at each cell $C=C_{i,j,k}$. A cell average is defined as the integral of the quantity over the cell, i.e., the cell average of the quantity $\textbf{U} = \textbf{U}_{i,j,k}$ is
\begin{equation}
    \overline{\textbf{U}} = \frac{1}{|C|}\int\limits_C \textbf{U} d\textbf{x},
\end{equation}
where $d\textbf{x} = dx\,dy\,dz$. 
By taking the integral form of Equation~\ref{eq:consLaw}, along with the concept of the cell average and the divergence theorem, we have
\begin{equation}
    \frac{\partial \overline{{\bf U}}}{\partial t} + \int\limits_{\partial C}{\bf F}\cdot{\bf n}_d\,d{\bf x} ={\bf 0},
\end{equation}
where ${\bf n}={\bf n}_d$ is the normal vector to the interface of the cell in direction $d$ and assuming that the boundary terms vanish. The integral on the left-hand side of the above equation means that the flux ${\bf F}$ must be evaluated on the interface of cell $C$, denoted by $\partial C$, projected onto the normal vector. The flux function can be numerically computed for each cell center, but it is not possible to analytically estimate it on the interface of the cells. 
In this context, it is necessary to make approximations of the flux. Here we have chosen the Harten-Lax-van Leer-Discontinuities (HLLD) Riemann solver \cite{Kusano:2005}. 

\begin{figure}[H]
\centering
\begin{tikzpicture}
%\draw[step=0.2cm,color=green, very thin] (-4,-4) grid (4, 4);
% Draw a line at 30 degrees and of length 3
\draw [->, very thick](0,0) -- (0:3.5cm) node[right]{$x$};
\draw (0,0) -- (18:3cm) node[right]{$\lambda_{C_{f^+}}$};
\draw (0,0) -- (36:3cm) node[right]{$\lambda_{C_{a^+}}$};
\draw (0,0) -- (54:3cm) node[right]{$\lambda_{C_{s^+}}$};
\draw (0,0) -- (72:3cm) node[right]{$\lambda_{C_{e}}$};
\draw [->,very thick](0,0) -- (90:3cm) node[above]{$t$};
\draw (0,0) -- (112.5:3cm) node[right]{$\lambda_{C_{s^-}}$};
\draw (0,0) -- (135:3cm) node[right]{$\lambda_{C_{a^-}}$};
\draw (0,0) -- (157.5:3cm) node[right]{$\lambda_{C_{f^-}}$};
\draw [very thick](0,0) -- (180:3cm);
\end{tikzpicture}
\caption{Illustration of the MHD Riemann fan.}
\label{fig:RF}
\end{figure}
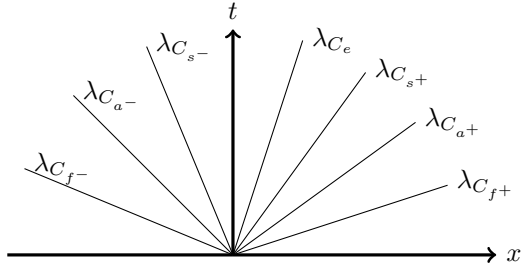

The HLLD scheme is useful to solve isolated discontinuities in the solution of the MHD system and it preserves the positivity of the quantities. 
The MHD eigensystem has seven eigenvalues, related to the fast magnetoacoustic waves ($C_f$), Alfv\`en waves ($C_a$), slow magnetoacoustic waves ($C_s$), and entropy ($C_e$). These waves compose the MHD space-time Riemann fan, illustrated in Figure~\ref{fig:RF}, in which they are organized according to their value regarding to their origin. 
The HLLD flux is efficient to solve isolated discontinuities and constructed by the means of an approximate Riemann problem in a four-state Riemann fan, divided by two Alfv\`en and one entropy wave.
This type of numerical flux is computed by sweeping the cells in each direction individually, a technique called splitting. 
Thus, the flux in $x$-direction, for instance, does not depend on the fluxes in the $y$ and $z$ directions. A brief summary of this Riemann solver can be found in \cite{Gomesetal:2015}. 
To calculate the numerical flux, the conservative variables are needed and we reconstruct them with the Monotonized Central (MC) limiter, which improves the accuracy of the solution by achieving a second order approximation.

\subsection{Adaptive multiresolution analysis}
After performing the discretization and flux approximations, we use multiresolution analysis to adapt the computational mesh to the numerical solution of the model. 
The multiresolution analysis comes from the idea that data can be decomposed into several levels of refinement $\ell$, from the coarser to the most refined level $L$, i. e., $0\leq \ell \leq L$. For details on multiresolution analysis we refer to \cite{Harten:1996}, or to textbooks \cite{Daubechies:1992,meyer1992wavelets, Mallat:1999}.
In a more particular way, $\overline{\bf U}^\ell$ represents a set of cell averages of the discrete data ${\bf U}$ at a refinement level $\ell$, with a total of $2^{3\ell}$ cells in three space dimensions.

When the data are composed into many levels, for instance, $\overline{\bf U}^{\ell-1}$, $\overline{\bf U}^\ell$, $\cdots$, $\overline{\bf U}^L$, it is possible to navigate between each level $\ell$ by defining prediction and projection operators such as
\begin{subequations}
    \begin{eqnarray}
    \mathcal{P}^{\ell\rightarrow \ell + 1}&:& \overline{\bf U}^\ell \longrightarrow \overline{\bf U}^{\ell+1}\\
    \mathcal{P}^{\ell +1\rightarrow \ell}&:& \overline{\bf U}^{\ell +1} \longrightarrow \overline{\bf U}^{\ell}.
    \end{eqnarray}
\end{subequations}
The projection operator $\mathcal{P}^{\ell\rightarrow \ell + 1}$ is exact and unique and it allows us to project the cell averages onto a coarser level. In contrast, the prediction operator performs an approximation of the cell averages in a more refined level, and it is not unique -- our choice here is linear interpolation. 

From the prediction process, we can compute the interpolation errors ${\bf d}^\ell={\bf d}^\ell_{i,j,k}$ at level $\ell$. These values are called details or wavelet coefficients and they have the information about the local regularity of the data. The cell averages together with the details at level $\ell +1$ enable us to obtain the cell averages at level $\ell + 1$, i. e., $\overline{\bf U}^{\ell+1}\leftrightarrow\{ {\bf d}^{\ell+1}, \overline{\bf U}^{\ell} \}$. 
More generally, it is possible to obtain the cell averages in the most refined level by the meaning of the cell averages at the coarsest level along with the details of every level between them
\begin{equation}
    L\leftrightarrow 0\,:\, \overline{\bf U}^{L}\leftrightarrow\{ {\bf d}^{L}, {\bf d}^{L-1},\cdots, {\bf d}^1,\overline{\bf U}^{0} \}.
\end{equation}
The relation above is called multiresolution transform and it ensures the possibility of approximating from data at any level $\ell$ by the knowledge of the details and cell averages of interest. 
For the three-dimensional case, $7$ wavelet coefficient are needed to perform the prediction operation of a cell. This procedure is revisited and presented in detail in \cite{DGRSESAIM:2011}.

The idea of the adaptive multiresolution analysis for cell averages \cite{Harten:1996} is to decrease the complexity of the multiresolution transform by locally adapting the computational mesh to the structures present in the data. 
This is made possible when we introduce a threshold operator which determines where to refine the mesh according to the detail information and a threshold parameter $\epsilon^\ell$. 
When the details are larger than $\epsilon^\ell$, the computational mesh needs to be more refined locally; otherwise the mesh can remain coarser. 
This methodology enables the computational mesh to be efficiently adapted to the numerical solution. 
In the context of the MHD model, where we have $8$ independent variables, the refinement criteria are applied to the maximum wavelet coefficient of these variables, so that the adapted mesh can be understood as a union of the meshes of each variable. 
In this work we are interested in the level-dependent threshold parameter $\epsilon^\ell=\epsilon^0$, which is modified according to the level of refinement and defined as following
\begin{equation}
\epsilon^\ell=\frac{\epsilon_0}{|\Omega|} 2^{3(\ell-L+1)}, \;\;0\leq \ell \leq L-1,
\label{eq:epsilon}
\end{equation}
where $|\Omega|$ is the volume of the domain and $\epsilon^0$ is the initial threshold parameter. 
Starting from $\epsilon^0$, the threshold value $\epsilon^\ell$ will increase for increasing $\ell$.
The number of cells on the finest mesh is defined as $3^{2L}$, where $L$ is the finest scale level. 

\subsection{Time evolution} 
To advance the solution in time from $t^n$ to $t^{n+1}$, we use a second-order compact Runge--Kutta scheme. 
The time steps are computed in terms of the Courant number $\nu$, the space steps, and the fast magnetoacoustic wave. 
Considering $\mathbb{E}$ as a time evolution operator, the multiresolution transform operator $\bf M$ and the threshold operator $\mathcal{T}$, we can summarize the entire process described in this section as
\begin{equation}
    \overline{\bf U}^{n+1} = {\bf M}^{-1}\circ\left\{ \mathcal{T}\circ\left[ {\bf M}\circ \left( \mathbb{E}\circ \overline{\bf U}^{n}\right)\right]\right\},
\end{equation}
where $\circ$ denotes the composition of functions.
%In both models the simulation is evaluated with a $2^{nd}$-order method. In CARMEN code, we have a compact Runge-Kutta, while in FLASH code a one-step Hancock is used.
The system is completed by choosing suitable initial and boundary conditions. 

\subsection{CARMEN and FLASH codes}
The CARMEN code described in \cite{RSTB03} is further developed in order to evaluate the numerical solution of the compressible Navier-Stokes including five more differential equations. An adaptive multiresolution algorithm for cell averages is used in combination with the finite volume discretization. 
The implementation of the two-dimensional version of the MHD equations, the Riemann solvers, parabolic-hyperbolic divergence cleaning and the variable reconstruction are already published in \cite{Gomes:2012:AnMuAd,Gomesetal:2015}. 
In this work, we have implemented the three-dimensional part of the MHD model, including the physical and numerical fluxes.

The FLASH code is an open-access parallel multiphysics multiscale simulation code provided by the Flash Center for Computational Science of the University of Chicago \cite{DubeyFlash:2014}. 
As an astrophysical code, it can simulate a great variety of problems. Here we are interested in the ideal MHD simulations in the context of the non-adaptive approach of the FLASH code, i. e., the full mesh finite volume simulation, to create a benchmark for our adaptive multiresolution results. 
We use version~$4.3$ of the FLASH code.

The time evolution in the FLASH code uses a one-step Hancock scheme and the $8$-wave algorithm \cite{Powell:1999} for divergence cleaning with a finite-volume cell-centered method. The variables are reconstructed with the MC limiter. 
As the numerical schemes are not the same, we chose the previous settings of the FLASH code to ensure that both codes are second-order accurate.

\section{Results}
\label{sec:results}
Herein we present and discuss the numerical results of the ideal MHD model for three-dimensional simulations using as initial condition one and two-dimensional Riemann problems. In the first study we compare the results with an exact solution and in the latter we compare the results with a benchmark solution obtained using the FLASH code at the same refinement level. 
The visualization of the solution is shown along with its mesh and the errors from the comparison between the results.

\subsection{One-dimensional Riemann initial condition}

In this case, the initial condition comes from the idea of a Riemann problem, which is basically characterized by a function that assumes a determined constant value in each half of the domain. 
At the interface between the two parts of the domain there is a discontinuity that must be represented accurately. 
The following setting, which has an exact solution\footnote{https://web.mathcces.rwth-aachen.de/mhdsolver/}, is used as a benchmark for our results. The initial condition is based on \cite{Li:2005} and given in Table~\ref{tab:R1D}.

\begin{table}[H]\centering
    \caption{One-dimensional Riemann initial condition}
\begin{small}
  \begin{tabular}{@{}cccccccccccccccc@{}}
  \toprule
    & $\rho$ && $p$ && $v_x$ && $v_y$ && $v_z$ && $B_x$ && $B_y$ && $B_z$\\
    \cmidrule{1-16}
     $x\leq 0$ & 1.08 && 0.95 && 1.2 && 0.01 && 0.5 && 2.0/$\sqrt{4\pi}$ && 3.6/$\sqrt{4\pi}$ && 2.0$\sqrt{2\pi}$ \\
     $x>0$ & 1.0 && 1.0 && 0.0 && 0.0 && 0.0 && 2.0/$\sqrt{4\pi}$ && 4.0/$\sqrt{4\pi}$ && 2.0$\sqrt{2\pi}$ \\
	\bottomrule
  \end{tabular}
  \end{small}
  \label{tab:R1D}  
\end{table}

We use the following  parameters in the simulations:  physical time $t=0.1$, the Courant number $\nu=0.3$, $\gamma=5/3$, $\alpha=0.4$, refinement level $L=8$, computational domain $[-0.5,0.5]^3$, and the threshold parameters $\epsilon^0 = 0.1$.

This initial condition can be extended to two and three-dimensions by tensor products, so that we have constant values of the variables in quadrants and sub-cubes, respectively. 
This test is useful to verify if the fluxes are being computed properly in every direction. 
Thus, we simulate this Riemann problem in each of the three directions and compare the results with the exact solution.

\begin{figure}[H]
    \centering
    \begin{tabular}{cc}
        $\rho$ & $p$ \\
        \includegraphics[width=0.45\linewidth]{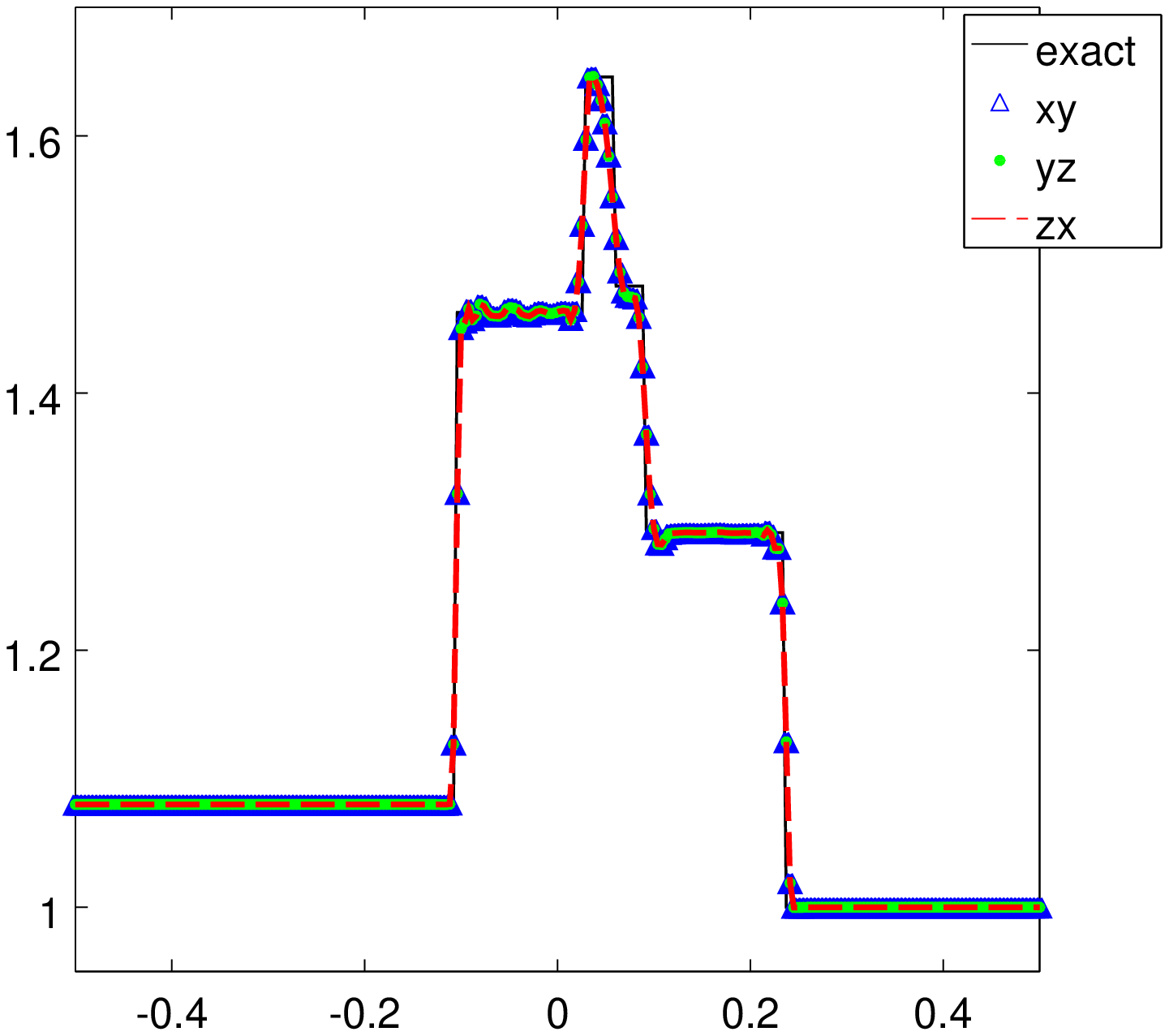} &             \includegraphics[width=0.45\linewidth]{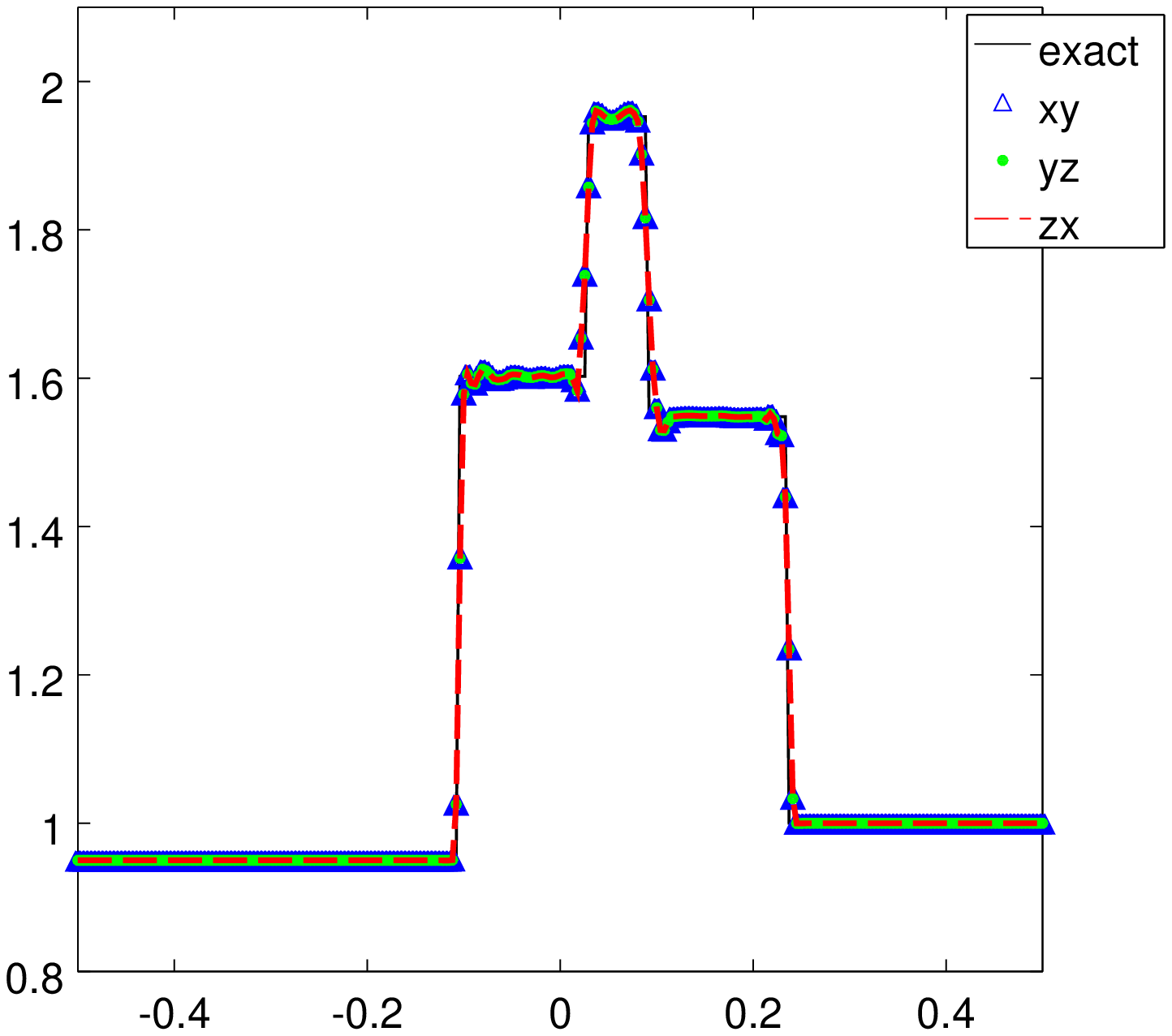}\\
        $u_x$ & $u_y$ \\
        \includegraphics[width=0.45\linewidth]{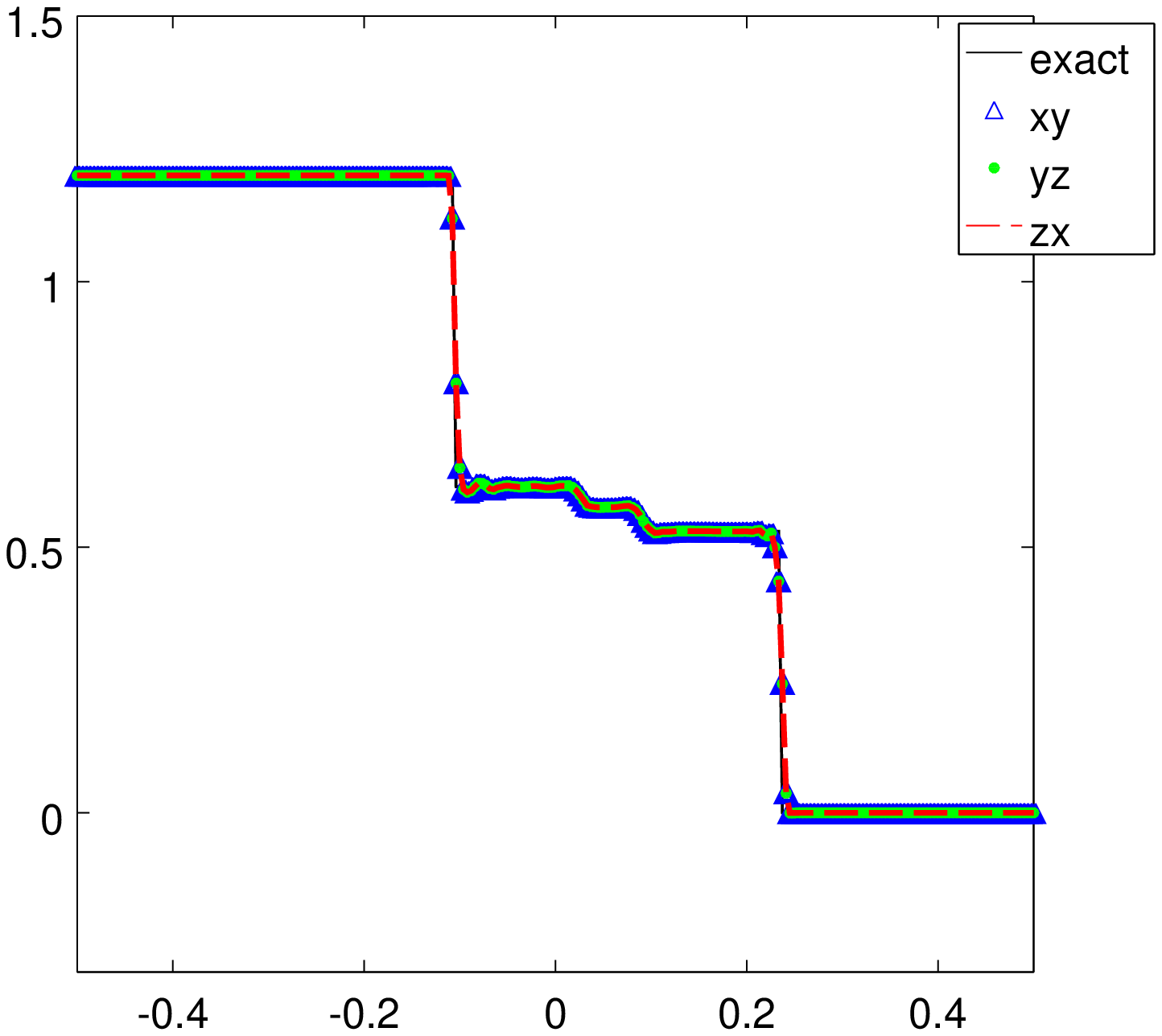} &             \includegraphics[width=0.45\linewidth]{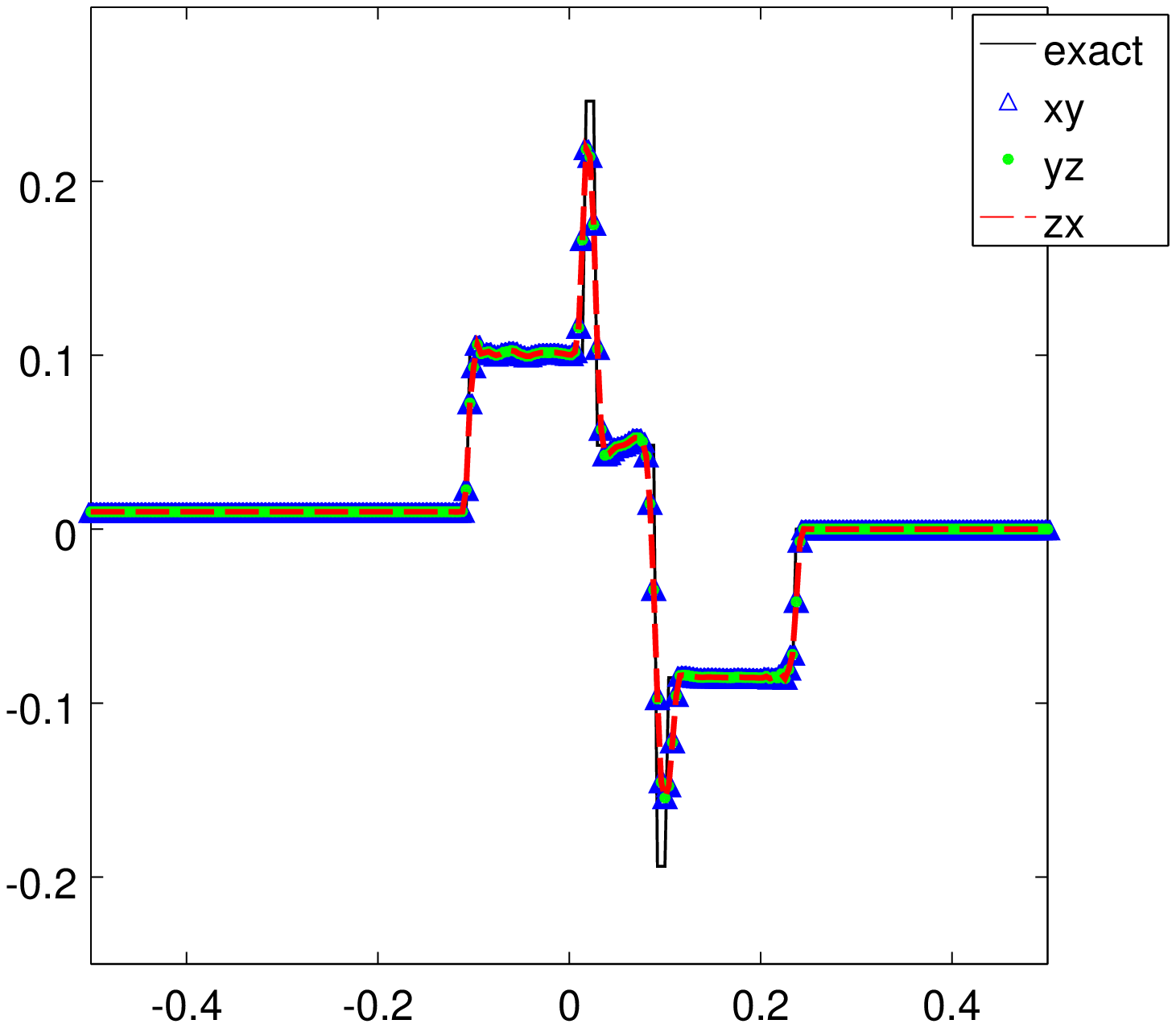}\\
        $B_y$ & $B_z$ \\
        \includegraphics[width=0.45\linewidth]{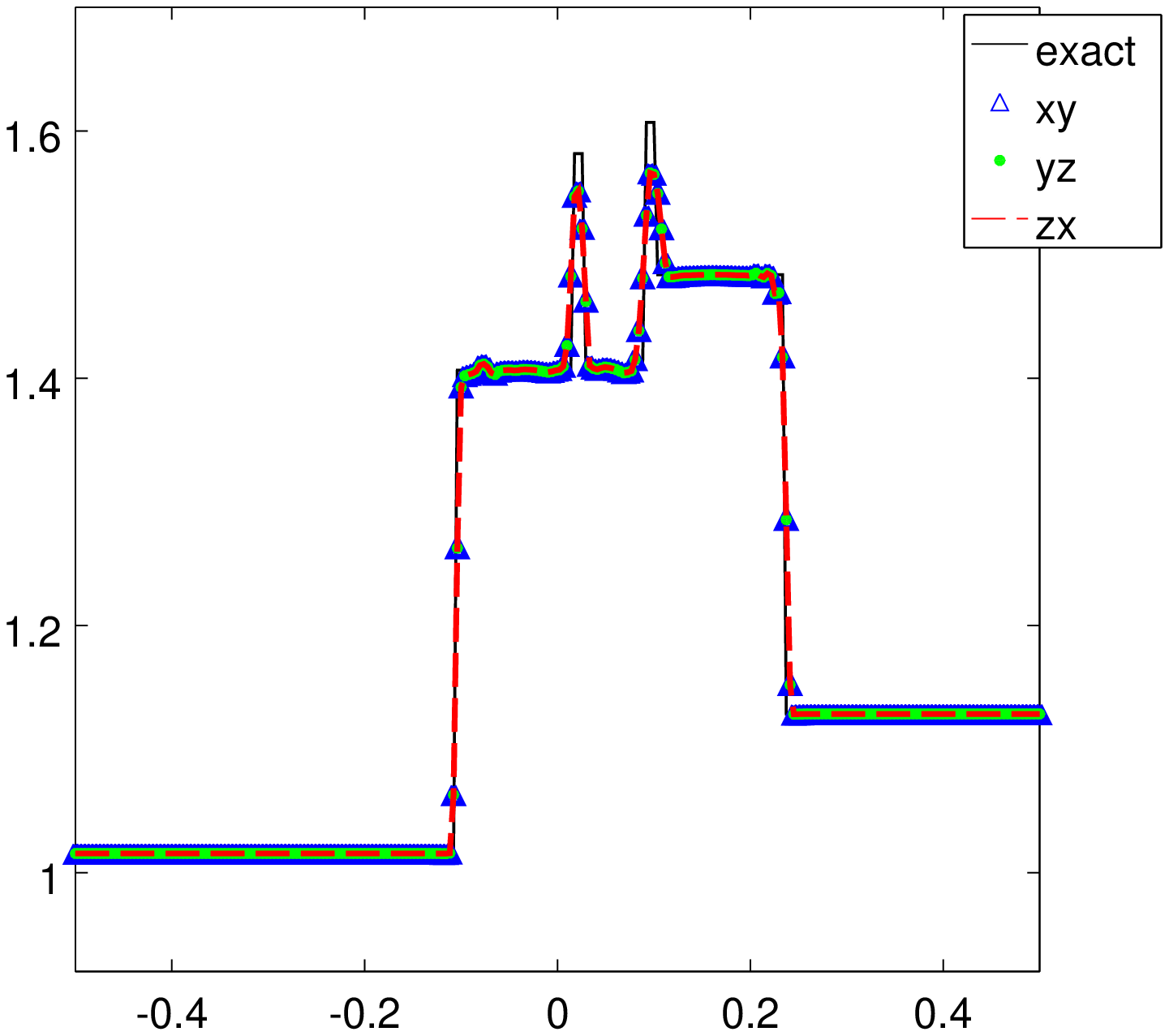} &             \includegraphics[width=0.45\linewidth]{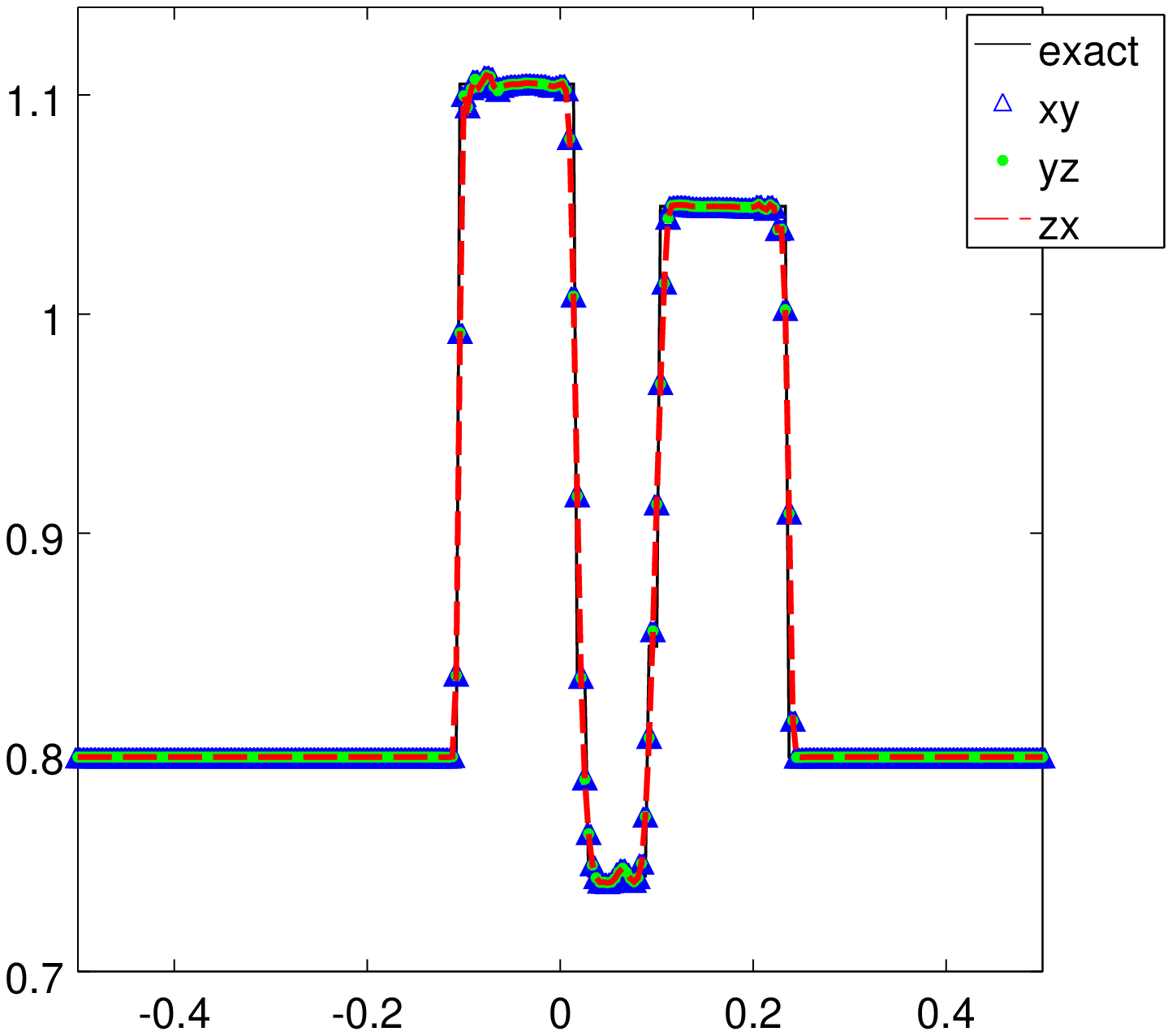}\\
    \end{tabular}
    \caption{Adaptive multiresolution computations. Shown are $\rho$, $p$, $u_x$, $u_y$, $B_y$, and $B_z$ obtained at $t=0.1$ and $L=8$, in the $xy$, $yz$, and $zx$ planes with $\epsilon^0 = 0.1$, superimposed is the exact solution (in black lines).}
    \label{fig:r1d3d}
\end{figure}

\begin{figure}
    \centering
    \begin{tabular}{ccc}
      \includegraphics[width=0.3\linewidth]{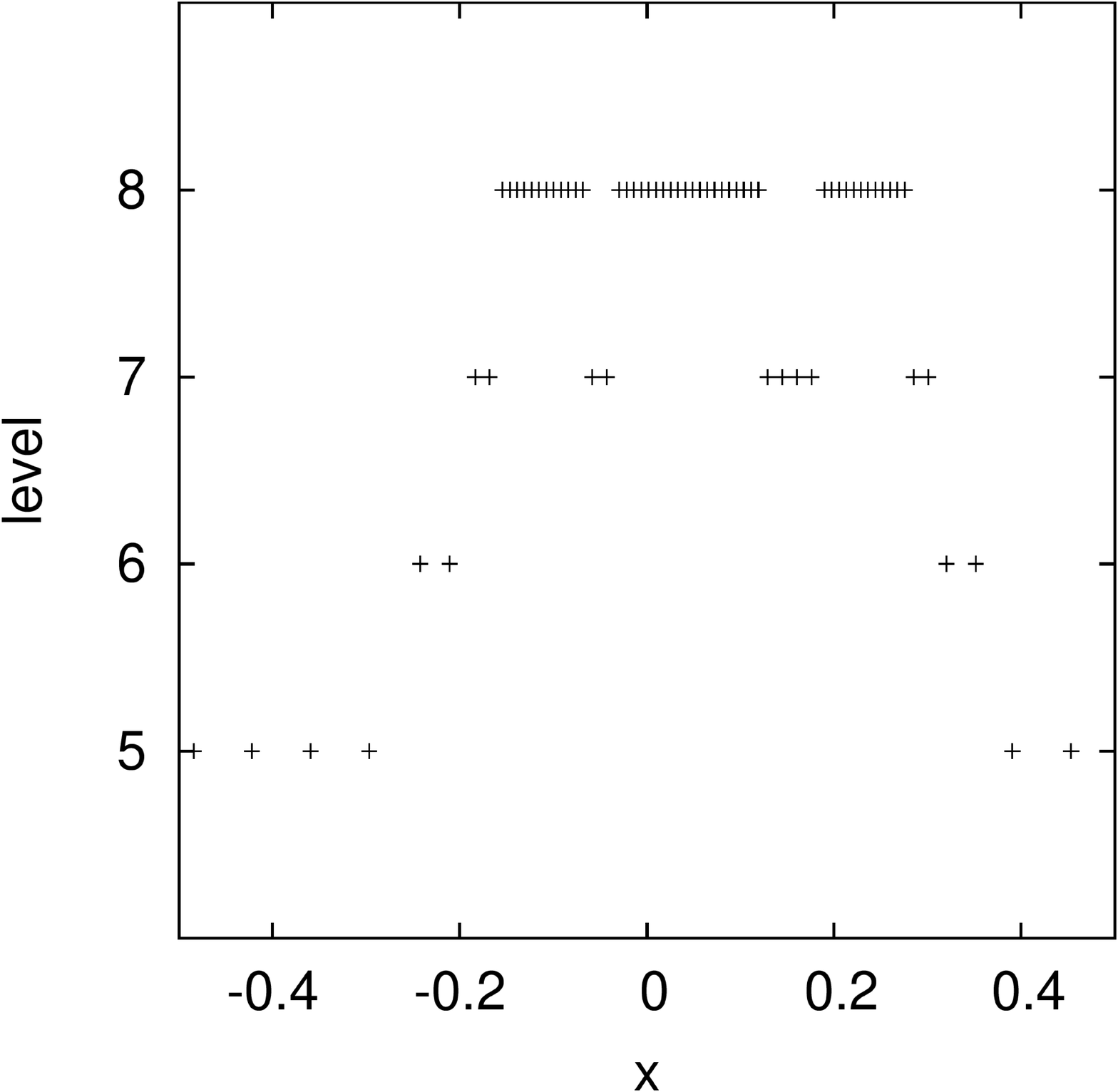} & \includegraphics[width=0.3\linewidth]{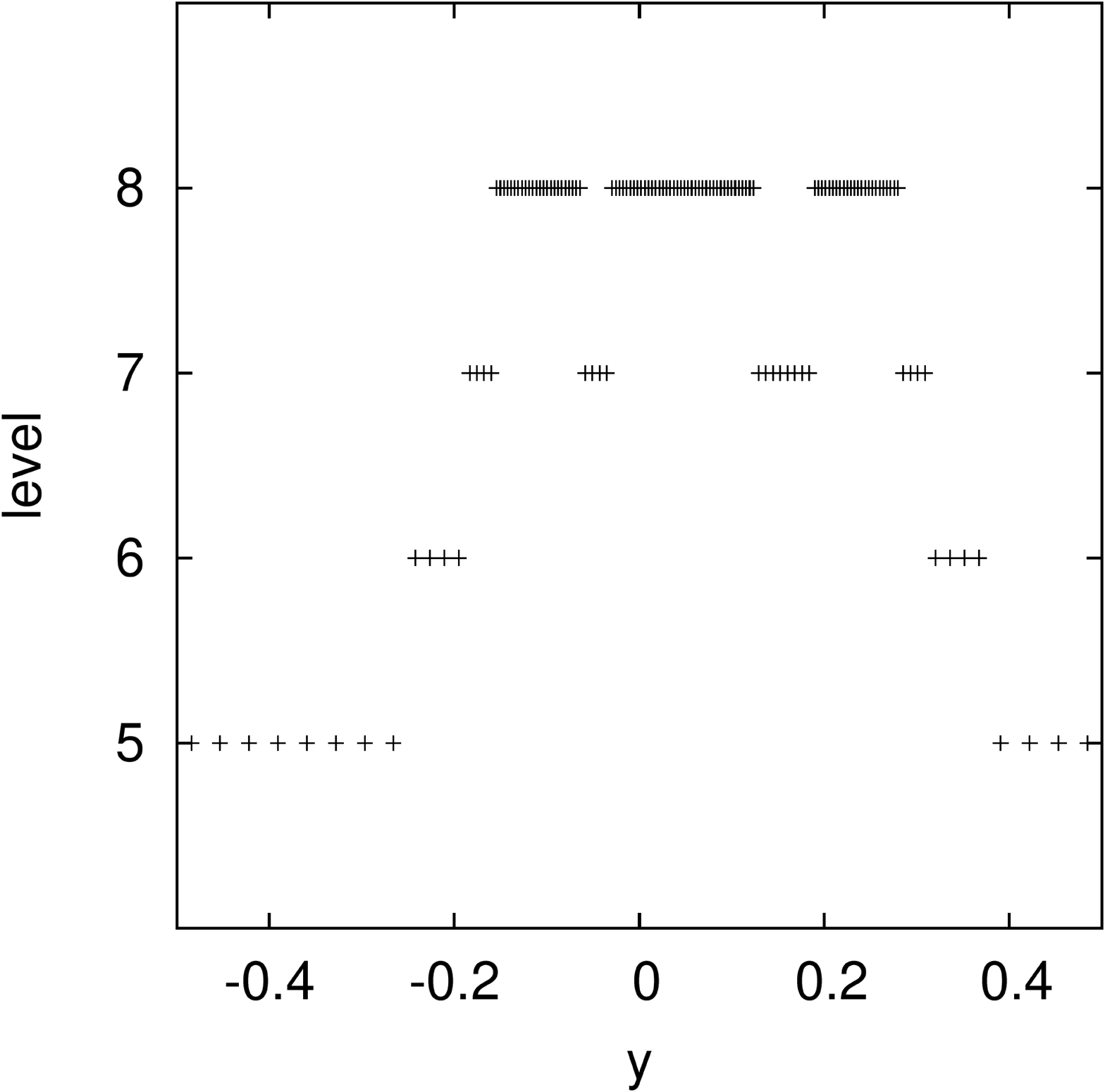} & \includegraphics[width=0.3\linewidth]{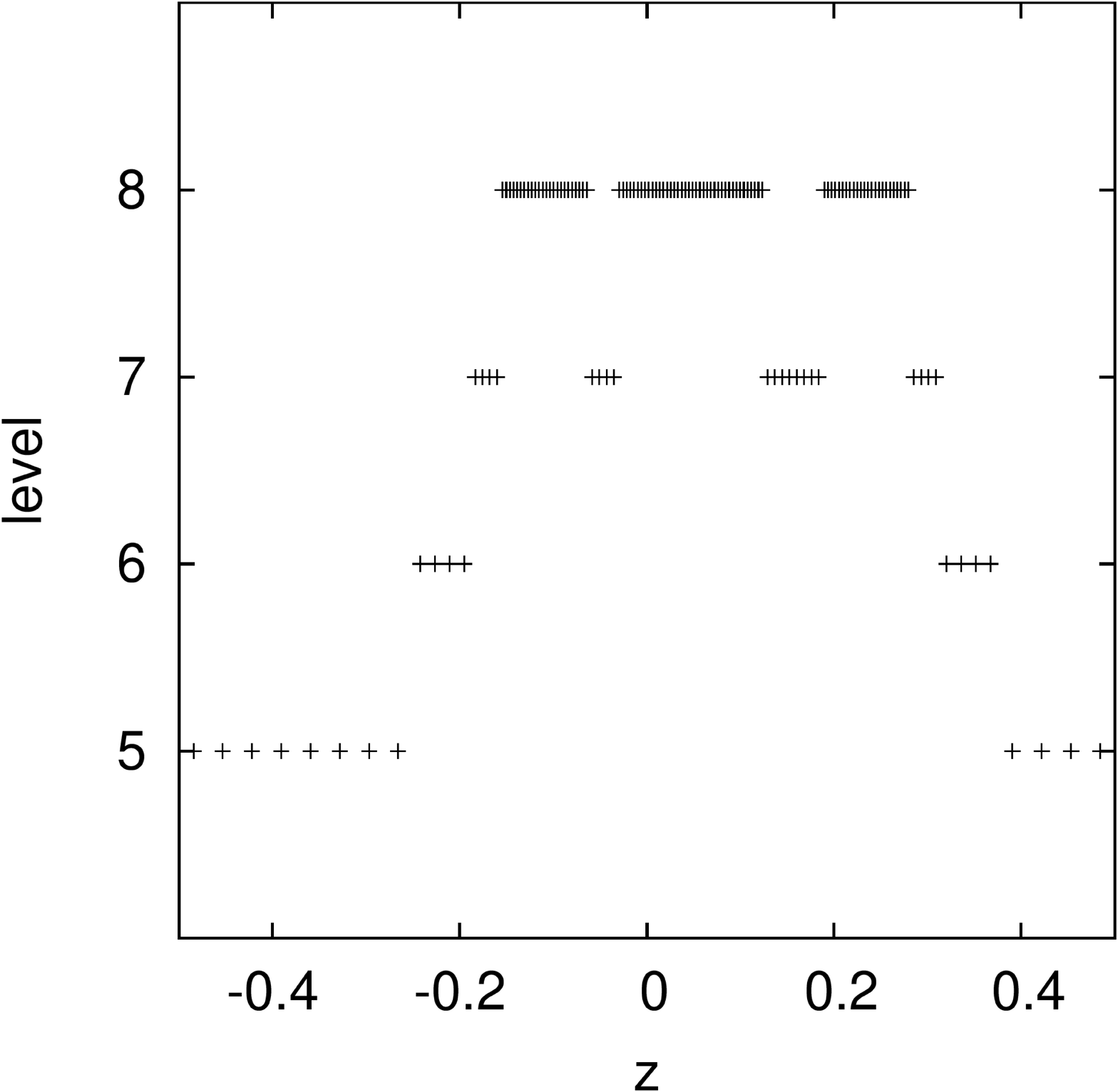}
    \end{tabular}
    \caption{Projection of the mesh in one-dimension in the $x$, $y$ and $z$ directions.}
    \label{fig:meshr1d3d}
\end{figure}

In Figure~\ref{fig:r1d3d} we can see the variables $\rho$, $p$, $u_x$, $u_y$, $B_y$, and $B_z$ of the solution of the Riemann problem obtained with the MR code for the threshold $\epsilon^0 = 0.1$ together with the exact solution at $t=0.1$ and for level $L=8$, which is computed in the planes $xy$, $yz$ and $zx$. 
The results show that the numerical solution converges to the exact solution for every variable of the model. 
The fluxes are being estimated in each direction of the domain and, from the figures, we can notice that their computations are reasonable. 
The one-dimensional projection of the mesh for this problem, composed by $36.6\%$ of the entire mesh, is shown in Figure~\ref{fig:meshr1d3d} for each case. The percentage of cells needed for the simulation over time is about $32\%$.

We can observe that the cells in the most refined level are limited to where the local structures are located. 
Levels $5$ and $6$, which are coarser, are representing the constant regions in the solution and more refined levels are needed only close to the localized structures. Table~\ref{tab:r1dNorm} shows the $\mathcal{L}^1$, $\mathcal{L}^2$ and $\mathcal{L}^\infty$ errors computed between the numerical and the exact solution for the one-dimensional Riemann problem for each variable of the MHD model.

\begin{table}[H]
    \centering
    \caption{Errors for the one-dimensional Riemann problem in three-dimensions for $L=8$ and $\epsilon^0=0.1$, obtained by comparing the results of the CARMEN code with the exact solution.}
    \begin{tabular}{@{}lccc@{}}
    \toprule
    Variables               &    $\mathcal{L}^1$ error     &    $\mathcal{L}^2$ error     & $\mathcal{L}^\infty$ error  \\ \cmidrule{1-4}
    $\rho$ & 4.5597$\cdot 10^{-3}$  &   1.0886$\cdot 10^{-3}$  &   1.4087$\cdot 10^{-1}$\\
    $p$    & 6.7027$\cdot 10^{-3}$  &   1.8227$\cdot 10^{-3}$  &   2.4520$\cdot 10^{-1}$\\
    $u_x$  & 3.5778$\cdot 10^{-3}$  &   1.3390$\cdot 10^{-3}$  &   2.4291$\cdot 10^{-1}$\\
    $u_y$  & 3.4091$\cdot 10^{-3}$  &   8.2139$\cdot 10^{-4}$  &   9.6215$\cdot 10^{-2}$\\
    $u_z$  & 2.5392$\cdot 10^{-3}$  &   6.6540$\cdot 10^{-4}$  &   7.8862$\cdot 10^{-2}$\\
    $B_x$  & 2.8865$\cdot 10^{-8}$  &   1.8041$\cdot 10^{-9}$  &   2.8865$\cdot 10^{-8}$\\
    $B_y$  & 4.7422$\cdot 10^{-3}$  &   1.1685$\cdot 10^{-3}$  &   1.5728$\cdot 10^{-1}$\\
    $B_z$  & 3.7892$\cdot 10^{-3}$  &   9.8471$\cdot 10^{-4}$  &   1.1345$\cdot 10^{-1}$\\
    \bottomrule
    \end{tabular}
    \label{tab:r1dNorm}
\end{table}

\subsection{Two-dimensional Riemann initial condition}
The two-dimensional Riemann initial condition follows the same idea as the one-dimensional Riemann problem. 
This time the domain is divided into four parts with constant values, analogous to a Cartesian plane. 
The original problem is in two-dimensions and based on the example presented in \cite{Dedneretal:2002}.
However, in this work we extend the problem to three-dimensions presented in Table~\ref{table:R2D}. 
%We do not have a exact solution for this case, thus in order to compare the obtained results we use the FLASH code as our benchmark.

\begin{table}[!htb]\centering
  \caption{Two-dimensional Riemann initial condition}
 \begin{small}
 \begin{tabular}{ccccccccc}
 \toprule
   Interval & $\rho$ & $p$ & $v_x$ & $v_y$ & $v_z$ & $B_x$ & $B_y$ & $B_z$\\
    \cmidrule{1-9}
       $x\geq 0$ and $y\geq 0$ & 1.0304 & 2.2874 & 1.4127 & -1.0146 & -1.0691 & 0.3501 & 0.5078 & 0.1576 \\
       $x< 0$ and $y> 0$ & 1.0000 & 2.4323 & 1.7500 & -1.0000 & 0.0000 & 0.5642 & 0.5078 & 0.2539 \\
       $x\leq 0$ and $y\leq 0$ & 1.8887 & 7.6110 & 0.1236 & -0.9224 & 0.0388 & 0.5642 & 0.9830 & 0.4915 \\
       $x>0$ and $y<0$ & 0.9308 & 2.1583 & 1.5639 & -0.4977 & 0.0618 & 0.3501 & 0.9830 & 0.3050 \\
	\bottomrule
  \end{tabular}
\end{small}
  \label{table:R2D}
  \end{table}
The simulation parameters are: the physical time $t=0.1$, the Courant number $\nu=0.3$, $\gamma=5/3$, $\alpha=0.4$, refinement level $L=8$, computational domain $[-1,1]^3$ and the threshold parameters $\epsilon^0=0.08$.

\begin{figure}
    \centering
    \begin{tabular}{cc}
        $p$ & Adapted mesh \\
        \includegraphics[width=0.45\linewidth]{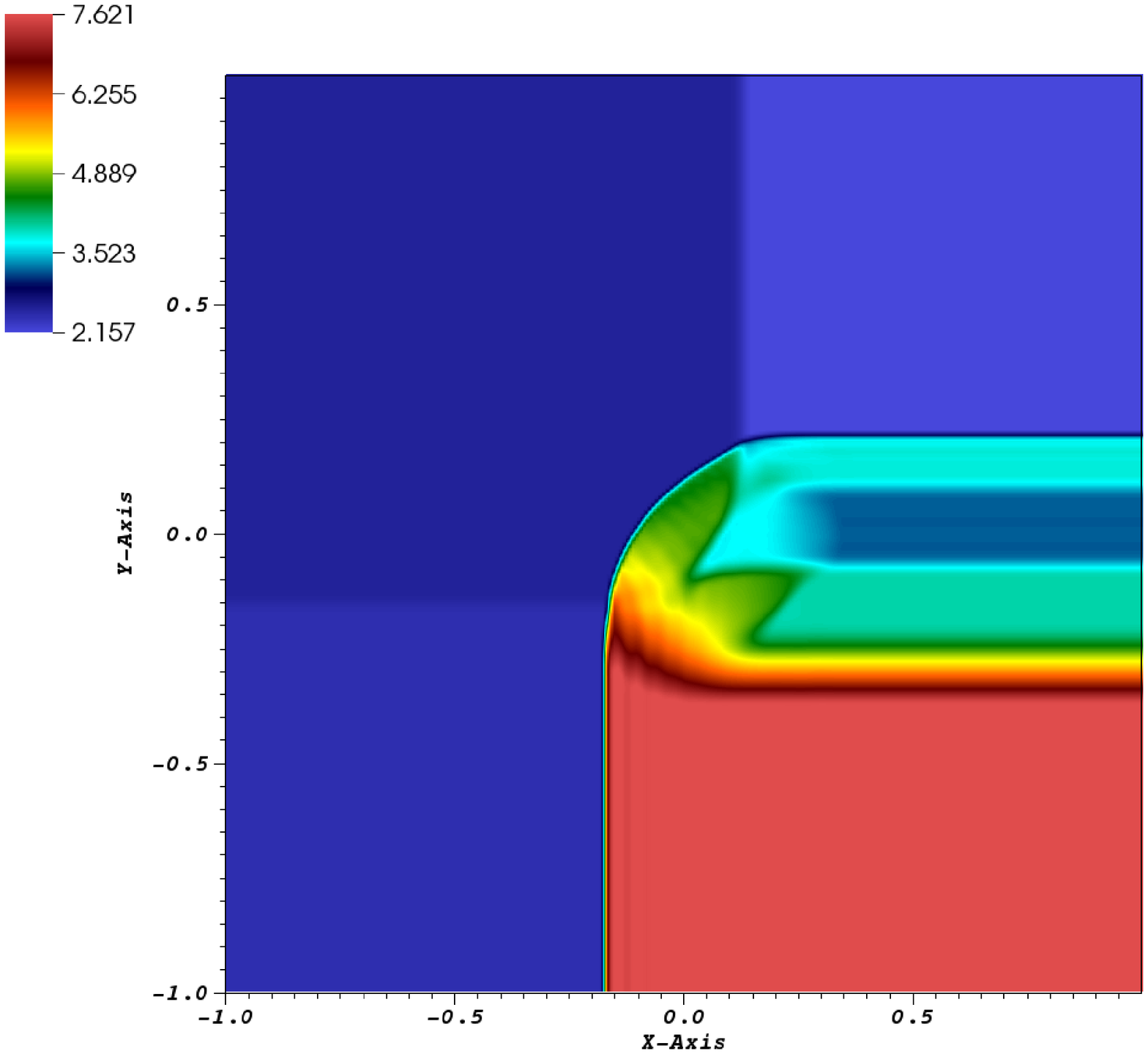} & 
        \includegraphics[width=0.47\linewidth]{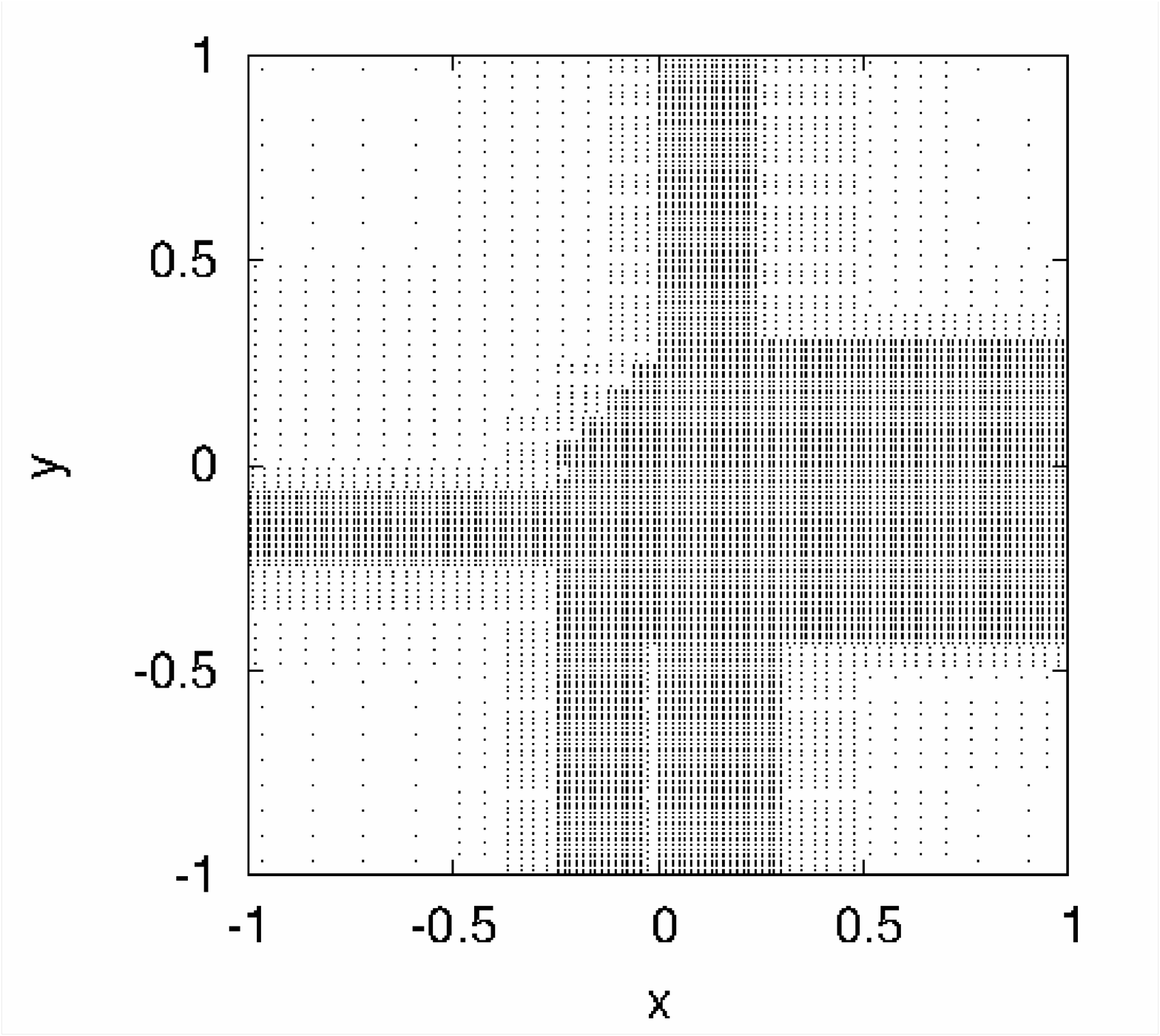} 
    \end{tabular} 
    \vspace{5mm}
    \caption{Variable $p$ and the projection of the mesh in two-dimensions ($x-y$ plane).}
    \label{fig:meshr2d3d}
\end{figure}

In Figure~\ref{fig:meshr2d3d} we show the variable $p$ and the corresponding adaptive mesh, composed by the cell-centers. The presented mesh is a projection of the 3D mesh onto the $xy$ plane and the darker regions mean the refinement is larger. The pressure is shown on Figure~\ref{fig:meshr2d3d} (left) and illustrates that the mesh is refined in regions of strong gradients. For the $yz$ and $zx$ plane simulation, the meshes are similar and at $t=0.1$ there are $40.8\%$ of the cells. 

In Table~\ref{table:R2D} the errors between the CARMEN and FLASH code results quantify the accuracy of each variable of the adaptive multiresolution simulation. The $\mathcal{L}^1$, $\mathcal{L}^2$, and $\mathcal{L}^\infty$ errors are approximately of the order of $10^{-3}$, $10^{-5}$, and $10^{-1}$, respectively. As expected, we have similar values for the $yz$ and $zx$ simulations.

\begin{table}[H]
    \centering
    \caption{Errors for the two-dimensional Riemann problem in three-dimensions for $L=8$ and $\epsilon^0=0.08$, obtained by comparing the results of CARMEN and FLASH codes.}
    \begin{tabular}{@{}cccc@{}}
    \toprule
    Variables               &    $\mathcal{L}^1$ error     &    $\mathcal{L}^2$ error     & $\mathcal{L}^\infty$ error  \\ \cmidrule{1-4}
  $\rho$ & 1.8604$\cdot 10^{-3}$ &   2.9467$\cdot 10^{-5}$ &   1.1195$\cdot 10^{-1}$  \\
  $p$    & 9.4471$\cdot 10^{-3}$ &   1.6666$\cdot 10^{-4}$ &   6.6384$\cdot 10^{-1}$  \\
  $u_x$  & 2.5411$\cdot 10^{-3}$ &   6.3577$\cdot 10^{-5}$ &   2.2006$\cdot 10^{-1}$  \\
  $u_y$  & 1.7706$\cdot 10^{-3}$ &   2.5164$\cdot 10^{-5}$ &   7.4709$\cdot 10^{-2}$  \\
  $u_z$  & 7.3367$\cdot 10^{-4}$ &   2.2670$\cdot 10^{-5}$ &   9.9015$\cdot 10^{-2}$  \\
  $B_x$  & 1.5709$\cdot 10^{-3}$ &   4.9229$\cdot 10^{-5}$ &   2.1496$\cdot 10^{-1}$  \\
  $B_y$  & 6.9340$\cdot 10^{-4}$ &   1.4988$\cdot 10^{-5}$ &   6.8804$\cdot 10^{-2}$  \\
  $B_z$  & 1.0497$\cdot 10^{-3}$ &   2.3591$\cdot 10^{-5}$ &   9.6685$\cdot 10^{-2}$  \\
    \bottomrule
    \end{tabular}
    \label{tab:kh-fv}
\end{table}

\begin{figure}[H]
    \centering
    \begin{tabular}{cc}
        $FLASH$ & $CARMEN$ \\
        \includegraphics[width=0.38\linewidth]{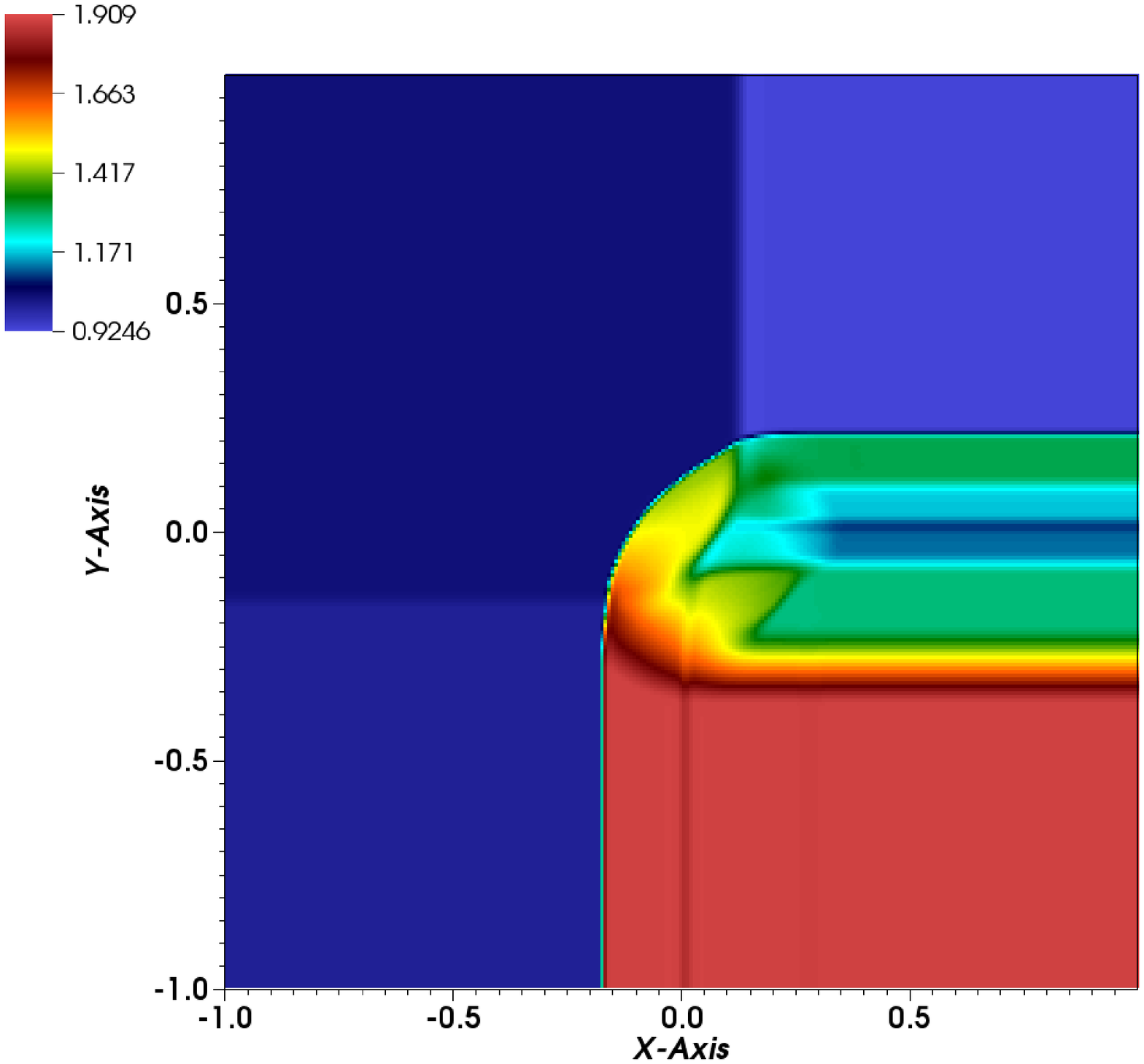} & 
        \includegraphics[width=0.38\linewidth]{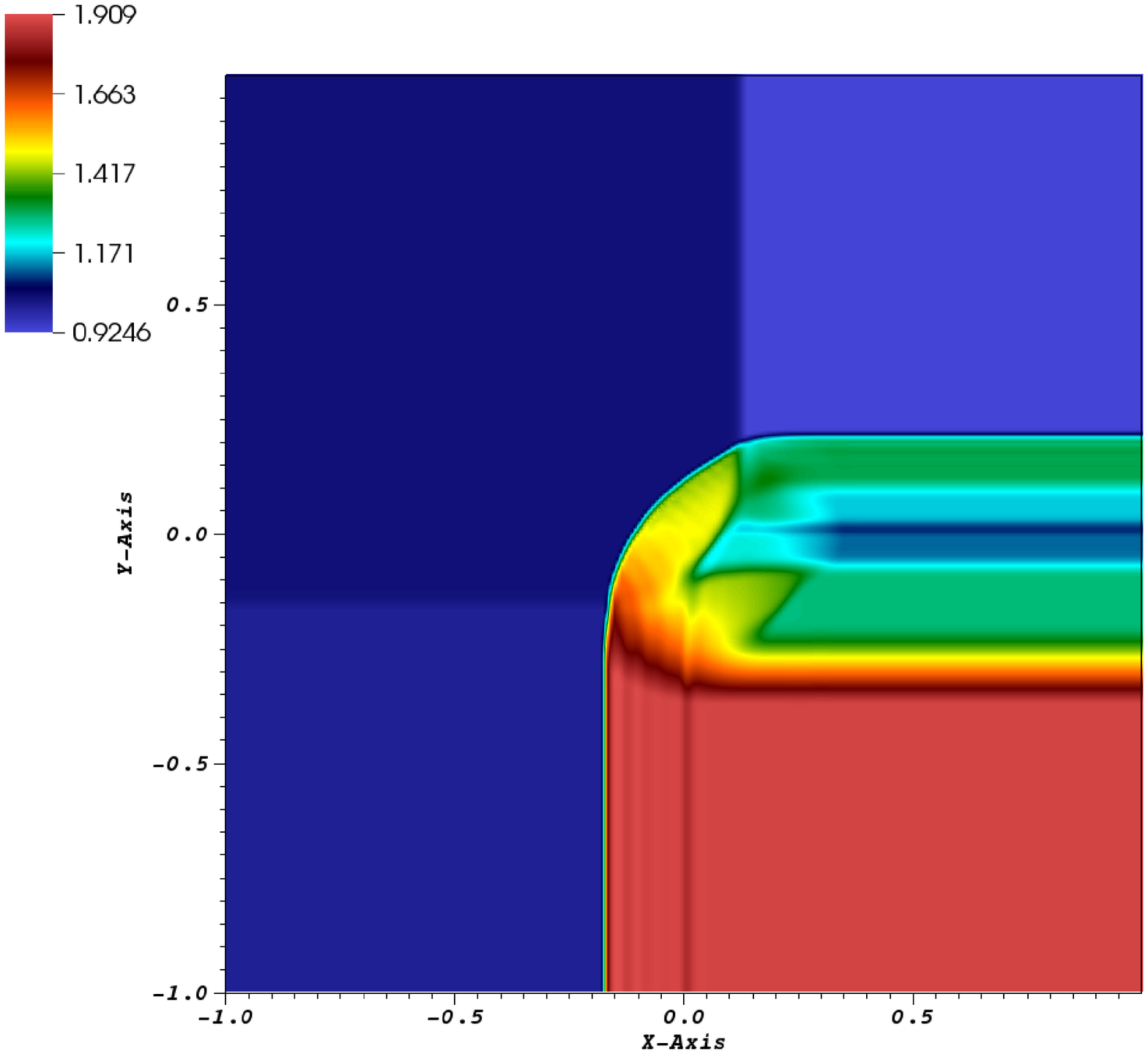}\\
		\includegraphics[width=0.38\linewidth]{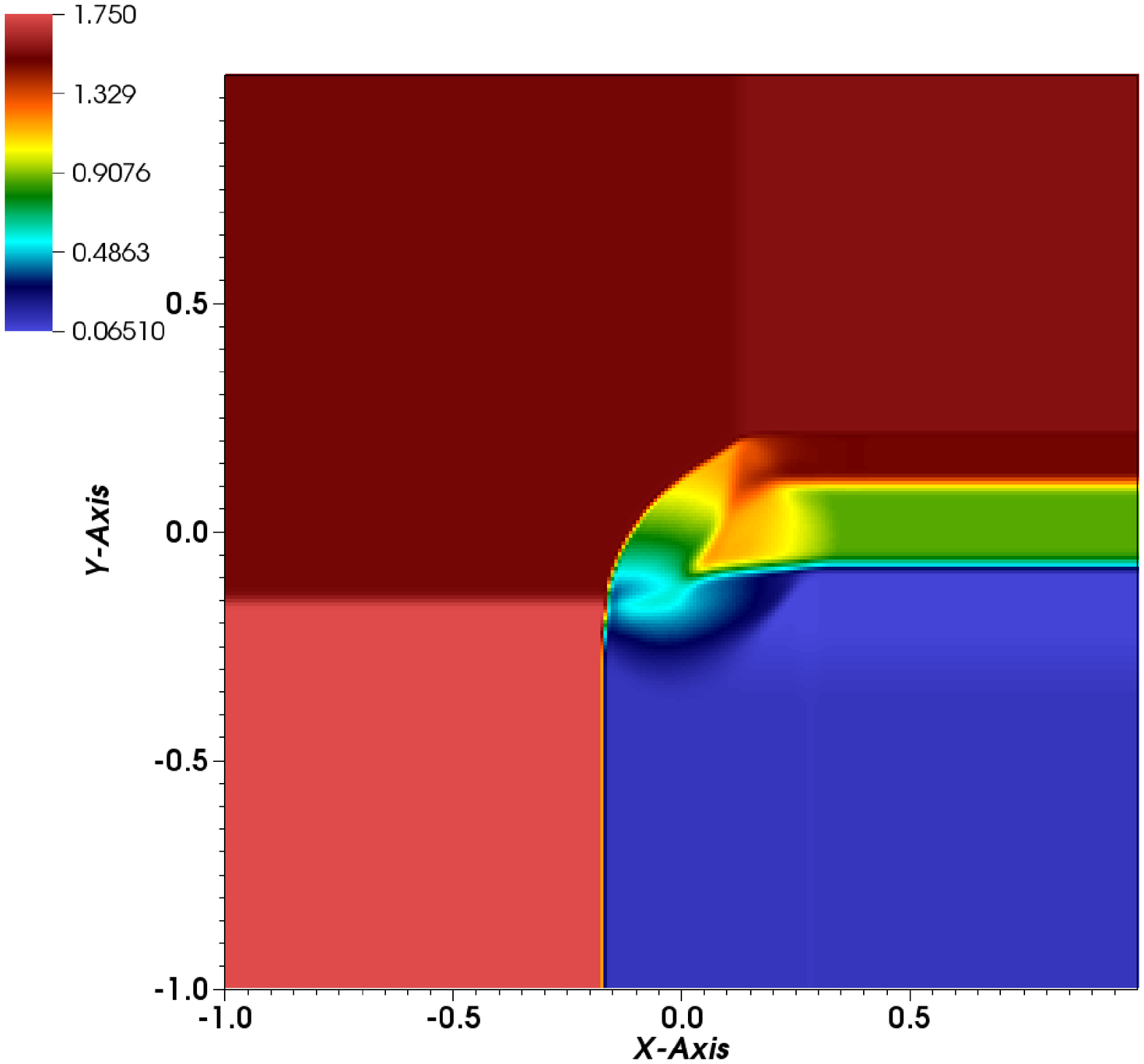} & 
        \includegraphics[width=0.38\linewidth]{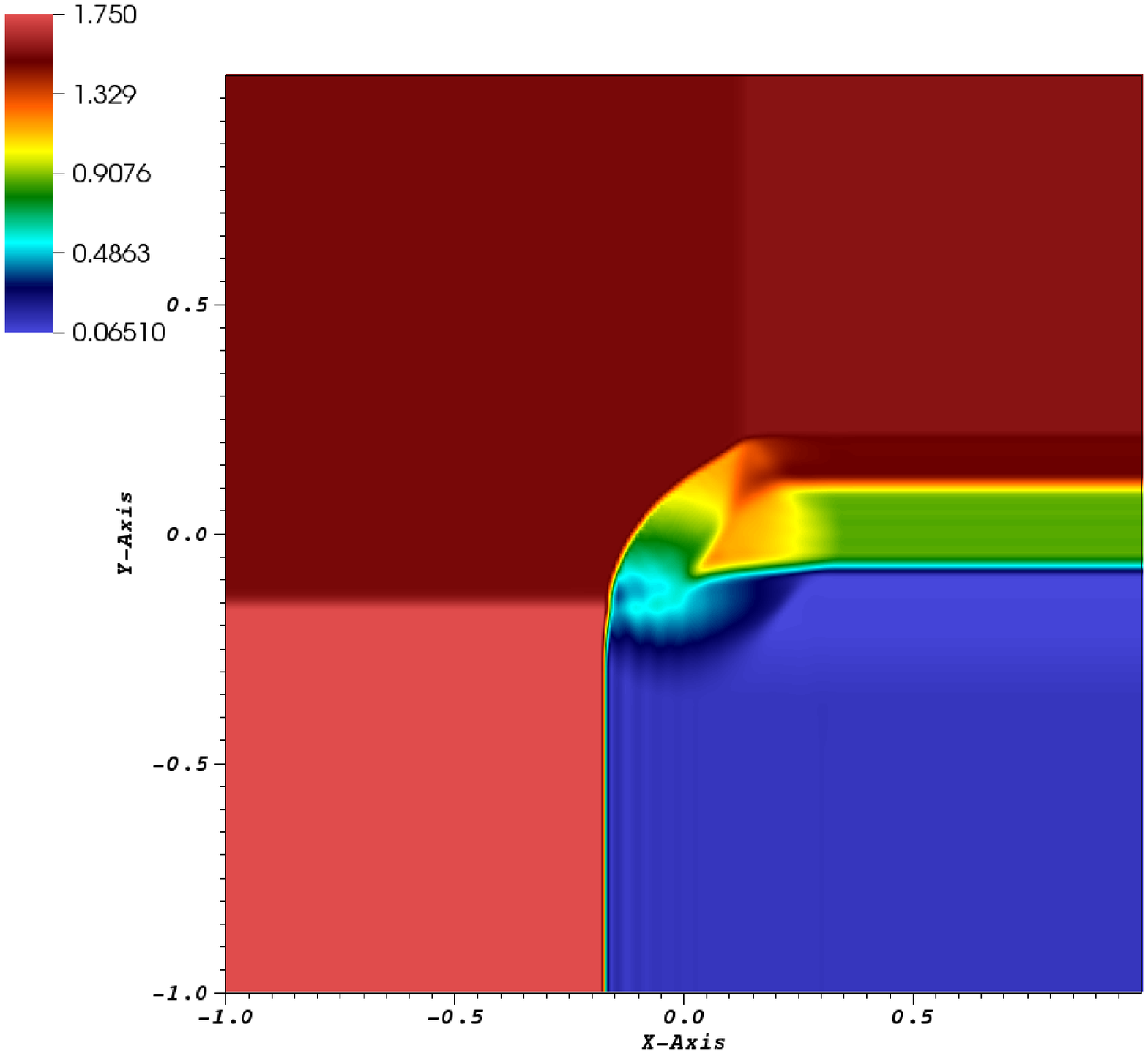}\\
        \includegraphics[width=0.38\linewidth]{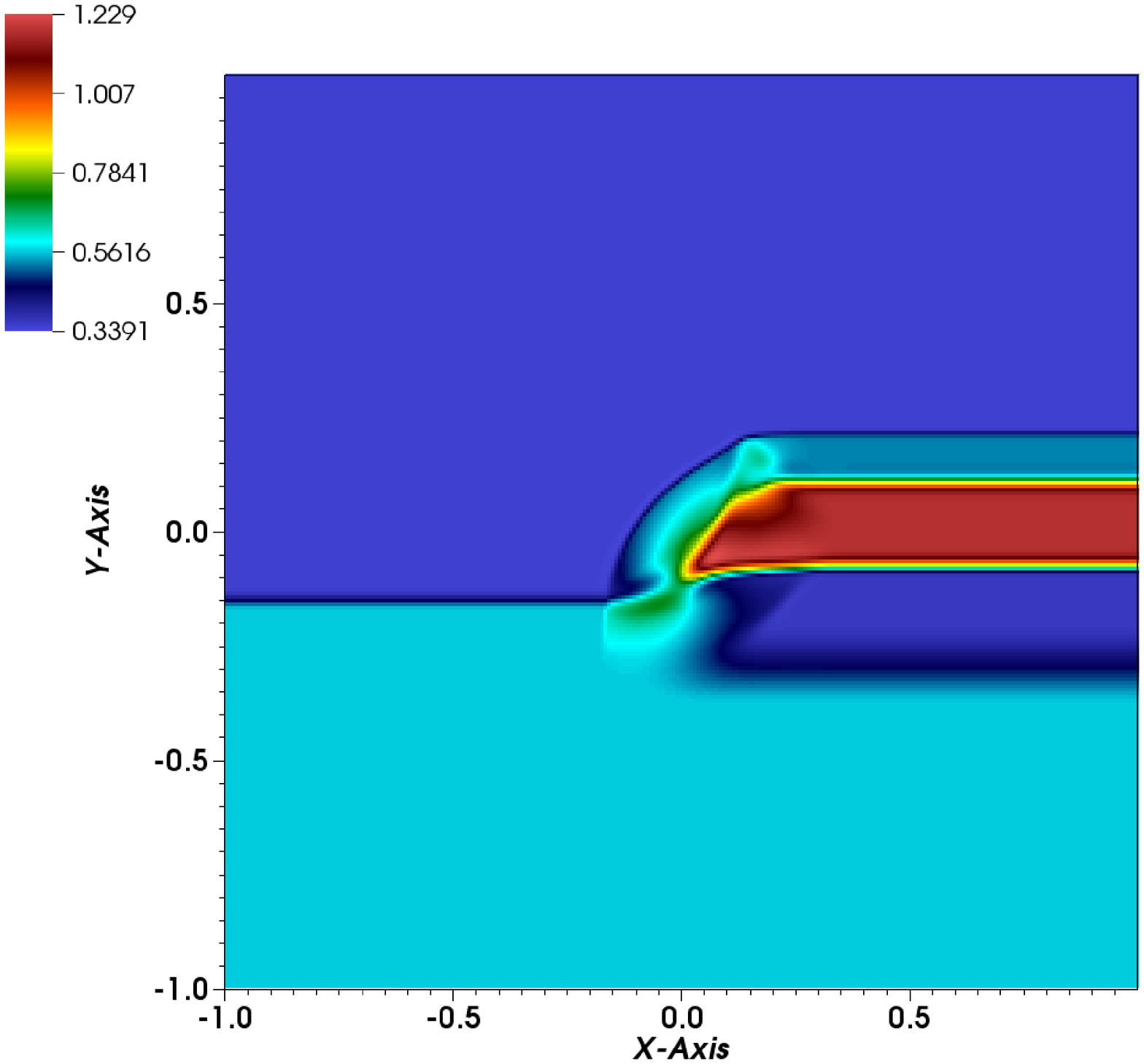} & 
        \includegraphics[width=0.38\linewidth]{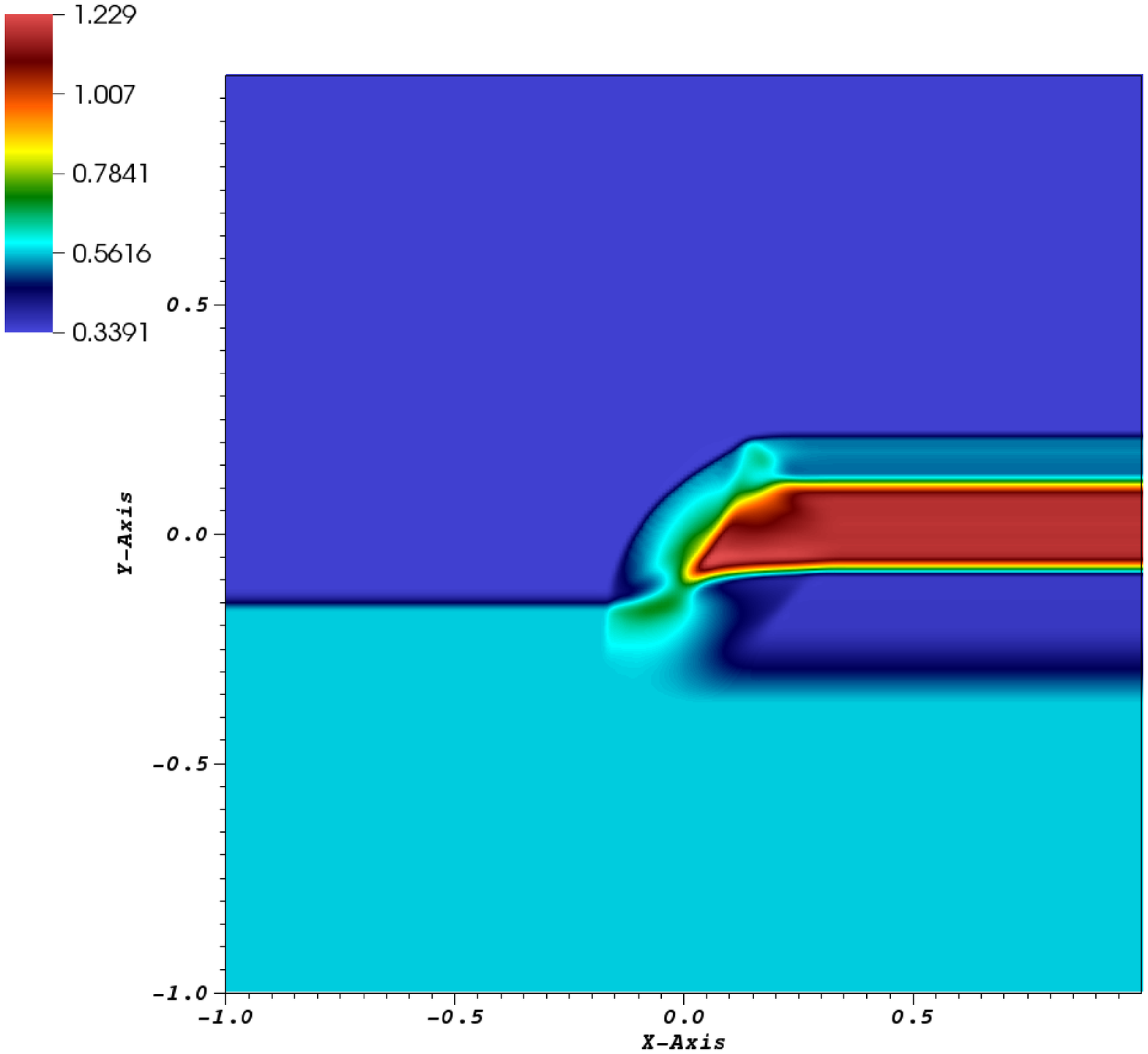}\\
    \end{tabular}
    \caption{Variables $\rho$, $u_x$, $B_y$, and $B_z$ obtained at $t=0.1$ and $L=8$, with the FLASH code and the CARMEN code for $\epsilon^0 = 0.08$ in the $xy$ plane.}
    %\caption{Variables $\rho$, $u_y$ and $B_x$ (from top to bottom) obtained at $t=0.5$ and $L=9$, with (a) FLASH code and (b) CARMEN code for $\epsilon = 0.25$.}
    \label{fig:r2d}
\end{figure}

The variables $\rho$, $v_x$, and $B_x$ are displayed in Figure~\ref{fig:r2d}, obtained with the FLASH (left) and CARMEN (right) codes. The solutions for both cases are similar using the same color scale. The structures located in the center of the domain are well represented and detached. The discontinuities between the constant regions are sharp in both directions in which they appear. For the problems computed on the $yz$ and $zx$ planes, the behavior of the solution is retained.

\section{Conclusions}
\noindent
We have presented three-dimensional simulations of the ideal MHD model with parabolic-hyperbolic divergence cleaning in the context of an adaptive multiresolution approach. 
It is shown that this type of approach indeed converges towards the expected results. 
With the MR method, the percentage of cells needed for the achievement of the solution is reduced to $32\%$ for the one-dimensional Riemann problem and about $39\%$ for the two-dimensional Riemann problem. 
Even with the reductions of cells in the simulation, the accuracy of the solution is well preserved.
The presented results are part of a verification of the 3D MHD CARMEN code and we can conclude that the numerical fluxes are computed properly in each direction and the MR approach is efficient to create an adaptive mesh for the studied problems.

\section*{Acknowledgments}
\noindent
The authors thank the Brazilian agencies CNPq ($306038/2015-3, 312246/2013-7$), FINEP/CT-INFRA ($01120527-00$), and FAPESP ($2015/50403-0, 2015/25624-2$)  projects for financial support. We thank Dr. Olivier Roussel for developing the original Carmen Code and for the fruitful scientific discussions. We also thank Eng. V. E. Menconi for his helpful computational assistance. 
The Flash software used in this work was in part developed by the DOE NNSA-ASC OASCR Flash Center
at the University of Chicago.
AKFG thanks CNPq for the PhD scholarship (project $141741/2013-9$).
MD thankfully acknowledges financial support from Ecole Centrale Marseille.
KS acknowledges financial support from the ANR-DFG, grant AIFIT, and support by the French Research Federation for Fusion Studies within the framework of the European Fusion Development Agreement (EFDA).

%\begin{thebibliography}{60}

\bibliographystyle{ieeetr} 
\bibliography{mhd}

\end{document}